\documentclass[12pt]{article}
 \usepackage{amsfonts}
 \usepackage{amssymb}
 \usepackage{amsmath,amsxtra}
\usepackage{amsthm}

\theoremstyle{plain}

\newtheorem{theorem}{Theorem}
\newtheorem{proposition}{Proposition}
\newtheorem{lemma}{Lemma}

\newtheorem{prop}{Proposition}

\theoremstyle{definition}
\newtheorem{deff}{Definition}

\newtheorem{rem}{Remark}

 \numberwithin{equation}{section}

\newcommand{\rt}{\rightarrow}

\newcommand{\mto}{\mapsto}
\newcommand{\bqa}{\begin{eqnarray}}
\newcommand\eqa {\end{eqnarray}}
\newcommand{\beq}{\begin{eqnarray}}
\newcommand{\beqn}{\begin{eqnarray}\nonumber}
\newcommand{\eeq}{\end{eqnarray}}
\newcommand{\be}{\begin{array}}
\newcommand{\ee}{\end{array}}
\newcommand{\bee}{\begin{equation}}
\newcommand{\eee}{\end{equation}}

 \newcommand{\pr}{\partial}
 \newcommand\na {\nabla}
 \newcommand{\p}{\partial}
 \newcommand{\Vol}{\operatorname{Vol}}

 \newcommand{\md}{\mathrm{d}}

\newcommand{\Id}{\mathrm{Id}}

 \newcommand{\Ber}{\mathrm{Ber}}

 \newcommand{\cA}{{\cal A}}
 \newcommand{\cD}{{\cal D}}
 \newcommand{\cM}{{\cal M}}

 \newcommand{\cI}{{\cal I}}
 \newcommand{\cX}{{\cal X}}
 \newcommand{\X}{{\cal X}}
 \newcommand{\cJ}{\mathcal{J}}
 \newcommand{\cL}{\mathcal{L}}

 \newcommand{\R}{{\mathbb R}}
 \newcommand{\Z}{{\mathbb Z}}

 \newcommand{\N}{{\mathbb N}}

   \newcommand{\g}{{\mathfrak g}}
   \newcommand{\gl}{{\mathfrak gl}}

   \newcommand{\sll}{{\mathfrak sl}}

\begin{document}

   \def\sp{\mathfrak sp}
   \def\sll{\mathfrak sl}
   \def\P{{\mathbb P}}
   \def\H{\mathbb H}
   \def\a{\alpha}
   \def\b{\beta}
   \def\t{\theta}
   \def\la{\lambda}
   \def\e{\epsilon}

 \def\gr{\g^{\scriptscriptstyle\mathrm{gr}}}
 \def\godd{\g_{\scriptscriptstyle 1}}
 \def\geven{\g_{\scriptscriptstyle 0}}
 \def\grodd{\gr_{\scriptscriptstyle 1}}
 \def\greven{\gr_{\scriptscriptstyle 0}}

  \def\sst{\scriptscriptstyle}
  \mathchardef\za="710B  
\mathchardef\zb="710C  
\mathchardef\zg="710D  
\mathchardef\zd="710E  
\mathchardef\zve="710F 
\mathchardef\zz="7110  
\mathchardef\zh="7111  
\mathchardef\zy="7112 

\mathchardef\zi="7113  
\mathchardef\zk="7114  
\mathchardef\zl="7115  
\mathchardef\zm="7116  
\mathchardef\zn="7117  
\mathchardef\zx="7118  
\mathchardef\zp="7119  
\mathchardef\zr="711A  
\mathchardef\zs="711B  
\mathchardef\zt="711C  
\mathchardef\zu="711D  
\mathchardef\zf="711E 
\mathchardef\zq="711F  
\mathchardef\zc="7120  
\mathchardef\zw="7121  
\mathchardef\ze="7122  
\mathchardef\zvy="7123  
\mathchardef\zvw="7124  
\mathchardef\zvr="7125 
\mathchardef\zvs="7126 
\mathchardef\zvf="7127  
\mathchardef\zG="7000  
\mathchardef\zD="7001  
\mathchardef\zY="7002  
\mathchardef\zL="7003  
\mathchardef\zX="7004  
\mathchardef\zP="7005  
\mathchardef\zS="7006  
\mathchardef\zU="7007  
\mathchardef\zF="7008  
\mathchardef\zW="700A  

\pagestyle{myheadings}\markright{Lie superalgebras of differential
operators}
\vskip 10mm

\title{Lie superalgebras of differential operators}

\author{Janusz Grabowski\thanks{The research of J. Grabowski
was supported by the Polish Ministry of Science and Higher Education
under the grant N N201 416839.}, Alexei Kotov, Norbert Poncin\thanks{
The research of A. Kotov and N. Poncin  was supported by UL-grant SGQnp2008.} }

\date{}

\maketitle

\abstract{\noindent We describe explicitly Lie superalgebra isomorphisms between the Lie superalgebras of first-order superdifferential operators on supermanifolds, showing in particular that any such isomorphism induces a diffeomorphism of the supermanifolds. We also prove that the group of automorphisms of such a Lie superalgebra is a semi-direct product
of the subgroup induced by the supermanifold diffeomorphisms and another subgroup which consists of automorphisms determined by even superdivergences. These superdivergences are proven to exist on any supermanifold and their local form is explicitly described as well.}

\vskip 4mm


 \vskip 2mm


\section{Introduction}

In \cite{GKP1} the authors have described all Lie superalgebra isomorphisms
$\cX (\cM_1)\to \cX (\cM_2)$ of the
Lie superalgebras of supervector fields on smooth supermanifolds $\cM_i$, $i=1,2$. It was shown
that every such an isomorphism is induced by a diffeomorphism of supermanifolds $\cM_1\to \cM_2$,
unless the dimension of the manifolds is $1|1$ or $0|2$. In the latter
low-dimensional cases there exist additional isomorphisms which have been also listed in \cite{GKP1}.

In the present paper the authors address a similar question, this time about isomorphisms of other natural and important Lie superalgebras growing on supermanifolds -- the Lie superalgebras $\cD^1(\cM)$ of first-order superdifferential operators. The structure of the Lie superalgebra automorphism group is much more complicated in this case, being a semi-direct product of the subgroup induced by the supermanifold diffeomorphisms and another subgroup which consists of automorphisms determined by even superdivergences.

\vskip 2mm \noindent The paper is organized as follows.
In section 2 we review some basic facts about (super)differential operators on a smooth supermanifold
and prove an important theorem (Theorem 1) about a representation of a supervector field
as a sum of supercommutators. We also explain why we focus on the first-order operators.

\vskip 2mm\noindent In section 3 we prove that,
if $\operatorname{dim}\cM$ differs from $0|1$, the algebra
$\cA(\cM)\subset\cD^1(\cM)$ of smooth (super)functions on $\cM$ can be characterized as the unique maximal super Lie ideal in $\cD^1(\cM)$ consisting of $\operatorname{ad}$-nilpotent elements.
In particular, any automorphism of the Lie superalgebra $\cD^1 (\cM)$ preserves $\cA(\cM)$,
thus induces an automorphism of the Lie superalgebra of supervector fields $\cX(\cM)\simeq\cD^1(\cM)/\cA(\cM)$.

\vskip 2mm\noindent In section 4 we prove that, according to the canonical splitting $\cD^1(\cM) =\cX(\cM)\oplus\cA(\cM)$, any automorphism of the Lie superalgebra $\cD^1(\cM)$ is a product
of an automorphism induced by a superdiffeomorphism of $\cM$ and
an automorphism of the form $(X+f)\mto (X+f+c(X))$,
determined by a 1-cocycle $c:\cX(\cM)\to\cA(\cM)$ on the Lie superalgebra
of vector fields $\cX (\cM)$, with values in the algebra of functions on $\cM$.
We also explain, why only those automorphisms of $\cX(\cM)$
which are induced by supermanifold diffeomorphisms can be extended from $\cX(\cM)$ to $\cD^1(\cM)$.

\vskip 2mm\noindent In section 5 we compute the 1st cohomology of the Lie superalgebra $\cX(\cM)$ of supervector fields with values in $\cA$. We show that every 1-cocycle is a combination of a closed super 1-form and
a fixed divergence.

\vskip 2mm \noindent
Finally, in section 6, we review necessary facts about measures (Berezinian volumes) and divergences on supermanifolds
and prove the existence of a divergence for each supermanifold.

\section{Sheafs of superdifferential operators}

Let us recall that a {\it smooth supermanifold} $\cM$ of dimension
$p|q$ is a (local) ringed superspace $(M,\cA)$ over a topological
space $M$ that is locally isomorphic to $(\R^p,C^{\infty\;p|q})$,
where, for any open subset ${\frak U}\,\subset\R^p$,
$C^{\infty\;p|q}({\frak U}\,):=C^{\infty}({\frak
U}\,)[\xi^1,\ldots,\xi^q]$ -- the $\xi^{\za}$ being formal
anticommuting generators. More precisely, we assume that $\cA$ is
a sheaf of associative supercommutative $\R$-algebras with unit.
The superalgebra $\cA(M)=\zG(M,\cA)$ of global sections of $\cA$
is the algebra $C^{\infty}(\cM)$ of {\it functions} of the
supermanifold $\cM$. It is well-known that, due to the local model
condition, the locality condition for the stalks is automatically
satisfied. Further, the considered data induce a smooth manifold
structure of dimension $p$ on $M$ and provide an embedding of the
classical manifold $M$ into the supermanifold $\cM$.

An important result of smooth supergeometry \cite{Gawedzki} asserts that there
exists a vector bundle $V$ of rank $q$ over $M$ such that $\cM$ is diffeomorphic as a supermanifold to $\Pi
V$, that is, to the total space of $V$ with the reversed parity of fibres. This implies that the algebra of
smooth functions on $\cM$ is isomorphic (as a commutative superalgebra) to the algebra of functions on $\Pi
V$, which is canonically identified with $\Gamma (\Lambda^\bullet V^*)$. This isomorphism is not canonical but
it gives us an identification
 \beq\label{isomorphism_functions}
  \cA=\cA_{\sst 0}\oplus\cA_{\sst 1}\simeq \Gamma (\Lambda^\bullet V^*)\,,
 \eeq
 with
 \beq\label{isomorphism_functions1}
  \cA_{\sst 0}=\bigoplus\limits_{i\ge 0} \cA^{\sst 2i}\,,\hspace{2mm}
  \cA_{\sst 1}=\bigoplus\limits_{i\ge 0} \cA^{\sst 2i+1}\,,\hspace{2mm}
  \mathrm{where}\hspace{1mm}\cA^{\sst k}\simeq \Gamma (\Lambda^{\sst k} V^*)\,.
 \eeq
In particular, this gives a (non-canonical) embedding of the algebra $C^\infty(M)$ into $C^\infty(\cM)$ and a super-version of the partition of unity.
\begin{rem} For simplicity, we assume in this text that $M$ is connected.\end{rem}

For any open subset $U\subset M$, we denote by
$(\operatorname{Der}\cA)(U)$ the $\cA(U)$-module
$\operatorname{Der}(\cA(U))$ of derivations of the superalgebra
$\cA(U)$. If $X\in (\operatorname{Der}\cA)(U)$, there is, in view
of the localization principle, for any open subset $V\subset U$, a
unique derivation $X|_V\in(\operatorname{Der}\cA)(V)$ such that
$(X f)|_V=X|_Vf|_V$, for all $f\in\cA(U)$. The assignment
$U\to(\operatorname{Der}\cA)(U)$ is actually a locally free sheaf
of $\cA$-modules, called the derivation sheaf
$\operatorname{Der}\cA$ of the structure sheaf $\cA$, or, also,
the tangent sheaf $T\cal M$ of the supermanifold $\cM$. The module
$(T\cM)(M)$ of global sections of the supervector bundle $T\cM$
is the $C^{\infty}(\cM)$-module ${\cal X}(\cM)$ of {\it vector
fields} of $\cM$ -- which carries an obvious Lie superalgebra
structure.\medskip

In the following we denote by $\underline{\operatorname{End}}(\cA(U))$ the $\cA(U)$-module of even and odd
$\R$-linear maps from $\cA(U)$ to itself. We can identify $\cA(U)$ with a subalgebra of
$\underline{\operatorname{End}}(\cA(U))$ using the left-regular representation $f\mapsto m_f$, $m_f(g)=fg$.
The $\cA(U)$-module of $k$-th order differential operators $\cD^k(U)$, $k\in\N$, is then defined inductively
by
\bee\label{Vin}\cD^k(U):=\{D\in
\underline{\operatorname{End}}(\cA(U)):[D,\cA(U)]\subset\cD^{k-1}(U)\},\eee
where $[-,-]$ is the supercommutator and where
$\cD^{-1}(U)=\{0\}.$

Of course, $\cD^0(U)=\cA(U)$, and thus 0-order operators are
local. This entails by induction that any superdifferential
operator is local. Indeed, if $D\in\cD^k(U)$, if the restriction
$f|_V$ of $f\in\cA(U)$ to an open $V\subset U$ vanishes, and if
$v\in V$, let $\zg\in \cA_0(U)$ be a bump superfunction with
support $\operatorname{supp}\zg\subset V$ (in the supercontext the
support can be defined as usual as the complement in $U$ of the
set of those points $u\in U$ for which the restriction of $\zg$ to
some neighborhood of $u$ vanishes) and restriction $\zg|_W=1$, for
some neighborhood $W\subset V$ of $v$, see localization principle
\cite[Corollary 3.1.8]{DAL}. It then follows from the defining
property of differential operators applied to $[D,\zg]f$, the
induction assumption, and the fact $\zg f=0$, that $(Df)|_W=0$. We
can now show that there exists, just as in the case of vector
fields, for any $D\in\cD^k(U)$ and any open $V\subset U$, a unique
$D|_V\in\cD^k(V)$ such that $(Df)|_V=D|_Vf|_V$, for all
$f\in\cA(U)$. Indeed, if $f\in\cA(V)$ and $v\in V$, it is possible
to choose a function $F\in\cA(U)$ (of the same parity as $f$) such
that $\operatorname{supp}F\subset\operatorname{supp}f$ and
$F|_W=f|_W$, for some neighborhood $W\subset V$ of $v.$ Locality
entails that $(DF)|_W\in\cA(W)$ and $(DF')|_{W'}\in\cA(W')$,
defined for two points $v,v'\in V$, depend only on $f$ and
coincide in the intersection $W\cap W'$. Thus these local
functions define a unique global function $D|_Vf\in\cA(V)$ such
that
$$(D|_Vf)|_W=(DF)|_W.$$ Since, obviously,
$D|_V\in\underline{\operatorname{End}}(\cA(V))$ (note that $D|_V$
has the same parity as $D$), it suffices -- to prove the above
claim -- to observe that, for any $f_1,\ldots,f_{k+1}\in\cA(V)$,
we have
$$[\ldots[[D|_V,f_1],f_2],\ldots,f_{k+1}]|_W=[\ldots[[D,F_1],F_2],\ldots,F_{k+1}]|_W=0,$$
with self-explaining notations.

In view of the just detailed restrictions of differential
operators, the assignment $U\to \cD^k(U)$ is a presheaf and
obviously also a sheaf -- as $\cA$ is a sheaf.

\begin{prop} For any $k\in\N$, the presheaf $\cD^k$ of $k$-th order superdifferential
operators over the base manifold $M$ (the body) of a smooth supermanifold
$\cM=(M,\cA)$ of dimension $p|q$ is a locally free sheaf of
$\cA$-modules, with local basis
$$\partial^{\za}_x\partial^{\zb}_{\xi}:=\partial_{x^1}^{\za^1}\ldots\partial_{x^p}^{\za^p}\partial_{\xi^q}^{\zb^q}\ldots\partial_{\xi^1}^{\zb^1},$$
where $(x,\xi)$ are local coordinates, $\zb_a\in\{0,1\}$, and
$|\za|+|\zb|\le k.$
\end{prop}

\begin{proof} The method used to prove local freeness of
the sheaf of vector fields goes through in the case of
differential operators. Let us give some details because of the
increased technicality.

If $M=U^{p|q}$ is a superdomain, if $D\in\cD^k(U)$ is of the type
\bee
\sum_{i=0}^kD^i=\sum_{i=0}^k\sum_{|\za|+|\zb|=i}D^i_{\za\zb}(x,\xi)\;\p_x^{\za}\p_{\xi}^{\zb
}\in\cD^k(U),\label{LocForm}\eee and if
$m_{\za\zb}=(1/\za!)x^{\za}\xi^{\zb},$ where the odd coordinates
are written in increasing order, then necessarily \bee
D^i_{\za\zb}=D^im_{\za\zb}=Dm_{\za\zb}-\sum_{j=0}^{i-1}D^jm_{\za\zb},\label{CoeffLocForm}\eee
and an induction on $i$ immediately shows that the coefficients
$D^i_{\za\zb}$, if they exist, are unique.

Take now an arbitrary $D\in\cD^k(U)$ and set
$\zD=D-\sum\in\cD^k(U)$, where $\sum$ denotes the {\small RHS} of
(\ref{LocForm}) with the coefficients defined in
(\ref{CoeffLocForm}). This operator $\zD$ vanishes by construction
on the polynomials of degree $\le k$ in $x,\xi$.

For any $f_1,\ldots,f_{\ell-1},h\in\cA(U)$, $\ell\ge k+1$, we have
\bee \zD(f_1\ldots f_{\ell-1}h)=\sum_{b=1}^{\ell-1}\sum\pm
f_{i_1}\ldots f_{i_b}\zD(f_{i_{b+1}}\ldots f_{i_{\ell-1}}h)+
F(h),\label{InducVanishPoly}\eee as immediately seen when
developing $F(h):=[\ldots[[\zD,f_1],f_2],\ldots,f_{\ell-1}]h.$ If
$\ell>k+1$, the term $F(h)$ vanishes, whereas in the case
$\ell=k+1$ it is given by $F(h)=F(1)h$. Equation
(\ref{InducVanishPoly}) shows that $\zD=0$ on any polynomial of
degree $k+1$, then, by induction, that $\zD=0$ on an arbitrary
polynomial in $x,\xi$. Further, this equation entails that
$\zD\cI^{k+c}_m\subset \cI^c_m$, $m\in U$, $c\ge 1$, where $\cI_m$
is the unique homogeneous maximal ideal of the stalk $\cA_m$.
However, in view of Hadamard's lemma, we can, for any $f\in\cA(U)$
and any $m\in U$, find a polynomial $P_{f,m}$ in $x,\xi$ such that
$f-P_{f,m}\in\cI^{k+q+1}_m$. It follows that $\zD
f=\zD(f-P_{f,m})\in\cI^{q+1}_m$, for all $m\in U$, so that $\zD
f=0.$
\end{proof}
\begin{rem}\label{order} It follows easily from the above proof that the order of a differential operator can
be determined locally by looking at supercommutators with the multiplications by coordinate functions:
$D\in\cD^k(U)$ if and only if, for a given system $(u^1,u^2,\ldots,u^{p+q})$ of local coordinates in $U$,
$$[[\ldots[[D|_U,u^{i_1}],u^{i_2}],\ldots],u^{i_{k+1}}]=0$$ for any sequence $i_1,\ldots,i_{k+1}\in\{
1,\ldots, p+q\}$.
\end{rem}
\begin{rem}
The idea of defining differential operators on an abstract (super)commutative (super)algebra by the formula
(\ref{Vin}) goes back to Grothendieck and Vinogradov \cite{Vi}.
\end{rem}
The super $\R$-vector space
$\underline{\operatorname{End}}(\cA(U))$ carries natural
associative and Lie superalgebra structures $\circ$ and $[-,-]$
(we often omit the symbol $\circ$). An induction on $k+\ell$
allows seeing that
$\cD^k(U)\circ\cD^{\ell}(U)\subset\cD^{k+\ell}(U)$ and
$[\cD^k(U),\cD^{\ell}(U)]\subset\cD^{k+\ell-1}(U)$, so that the
supervector space $\cD(U):=\cup_{k\in\N}\cD^k(U)$ of all
differential operators inherits associative and Lie superalgebra
structures that have weight $0$ and $-1$, respectively, with
respect to the filtration degree. It is easily checked that
$\cD:U\to\cD(U)$ (resp. $\cD^1:U\to\cD^1(U)$) is a locally free
sheaf of $\cA$-modules and associative and Lie superalgebras
(resp. of $\cA$-modules and sub Lie superalgebras) over $M$. The
algebra $\cD(M)$ (resp. $\cD^1(M)$) is the Lie superalgebra of
{\it differential operators} (resp. {\it first-order differential
operators}) of the supermanifold $\cM$. In the sequel we denote
this algebra also by $\cD(\cM)$ or even by $\cD$ (resp. by
$\cD^1{(\cM)}$ or $\cD^1$).\medskip

The usual splitting of the space of first-order differential
operators holds true in the supersetting.

\begin{prop} Let ${\cal M}=(M,{\cal A})$ be a smooth supermanifold. An endomorphism $D\in\underline{\mathrm{End}}({\cal A})$ is a first-order differential operator if and only if, for any
$f,g\in{{\cal A}}$, we have \bee\label{Char1stDiffOp} D(f g)=(D f) g+(-1)^{|D||f|} f (D g) -(D1)f g,\eee so
that supervector fields are those first-order differential operators $D$ that verify $D1=0$. Moreover, the
supervector space ${\cal D}^1$ admits a canonical splitting \bee {\cal D}^1={\cal A}\oplus {\cal X}\eee given
by $D\mapsto D1+(D-D1)$.\end{prop}

\begin{proof} To prove the direct implication (resp. converse implication), it suffices to compute $[[D,f],g](1)$ (resp. $[D,f](g)$), where $[-,-]$ denotes
the supercommutator of endomorphisms and where $1$ is the unit of
${\cal A}$. If $X\in{\cal X}$, we have $[X,f]=Xf\in{\cal
A}$, so that ${\cal X}\subset {\cal D}^1$. The second
claim now follows from Equation (\ref{Char1stDiffOp}). It entails
that the sum $\cal A+\cal X$ is direct. Finally, any $D\in{\cal
D}^1$ decomposes in the form $D=D1+(D-D1)$, where $D1\in {\cal
A}$ and $D-D1\in {\cal X}$.
\end{proof}
For purely even manifolds it is a well-known fact that the derived ideal $\cX(M)'=[\cX(M),\cX(M)]$ is the
whole Lie algebra $\cX(M)$ of vector fields of $M$ \cite{Grabowski}. More precisely, if $X\in\cX(M)$,
$\mathrm{supp}\,X\subset U$, $U$ open in $M$, then $X=\sum_{i=1}^n[X_i,Y_i]$, where $n$ is independent of the
considered $X$ and where $X_i,Y_i$ are vector fields of $M$ with support $\mathrm{supp}\,X_i$,
$\mathrm{supp}\,Y_i\subset U$ \cite{Pon04}. The next theorem extends these result to the supercontext and will
be used as a technical tool in the sequel.
\begin{theorem}\label{local} Let $\cM$ be a smooth supermanifold. Then, every $X\in\cX(\cM)$ can be written as a finite sum
of supercommutators,
\beq\label{L} X=\sum_i[X_i,Y_i]\,,\eeq where $X_i\in\cX(\cM)$
and $\,Y_i\in\cX_0(\cM)$. Moreover, if \ $\mathrm{supp}\,X\subset U$, $U$ open in $M$, then
$X_i,Y_i$ can be chosen so that $\mathrm{supp}\,Y_i\subset\mathrm{supp}\,X $ and $\mathrm{supp}\,(\zp X_i)\subset U$ for any $i$. In particular, the
derived algebra $\cX'(\cM)=[\cX(\cM),\cX(\cM)]$ equals $\cX(\cM)$.
\end{theorem}

\begin{proof} Using the theorem stating that any smooth supermanifold $\cM=(M,\cA)$ is
diffeomorphic to the supermanifold $\Pi V$ for some vector bundle $V$ over $M$, we can assume that $\cM=\Pi
V$. The algebras of functions and of vector fields of supermanifolds of the type $\zP V$ carry a
$\Z_2$-compatible $\N$-grading. We denote the parity-supergrading by subscripts and the $\N$-grading by
superscripts. The grading is recognized by the {\em Euler vector field} $\ze$, in the sense that $[\ze,X]=nX$
for $X\in\cX^n(\Pi V)$. Let us recall that supervector fields of degree 0 which, in local coordinates
$(x,\xi)$, read
$$X=\sum_iX^i(x)\p_{x^i}+\sum_{a,b}X^a_b(x)\xi^b\p_{\xi^a},$$
can be identified with the sections of the Atiyah algebroid of $V$. This identification is a Lie algebra
isomorphism, so that
$$0\to\ker\zp\to \cX^0(\zP V)\stackrel{\zp}{\to} \cX(M)\to 0$$ is a split short exact
sequence of Lie algebras, so that $\cX^0(\zP V)= \ker\zp\oplus\cX(M)$. Note that in coordinates
$$\pi\left(\sum_iX^i(x)\p_{x^i}+\sum_{a,b}X^a_b(x)\xi^b\p_{\xi^a}\right)=\sum_iX^i(x)\p_{x^i}$$ and the Euler
vector field reads $\sum_a\xi^a\p_{\xi^a}$.

\smallskip
We first prove the theorem for $X\in\cX(M)$, $\mathrm{supp}\,X\subset U$. According to \cite[Theorems (3.2)
and (4.1)]{Grabowski0}, $$X\in[\cX(M),[\cX(M),X]]\,,$$ so that there are $X_i,Y'_i\in\cX(M)$ such that
$$X=\sum_i[X_i,[Y'_i,X]]\,.$$
Putting $Y_i=[Y'_i,X]$, we can write $X$ in the form (\ref{L}) with
$\mathrm{supp}\,Y_i=[Y'_i,X]\subset\mathrm{supp}\,X$. Let us observe that this is also true for any
$X\in\cX(\Pi V)$. Indeed, if $X\in\cX^k(\Pi V)$, $k\ne 0$, then
$$X=\frac{1}{k}[\epsilon,X]\,,$$ where $\epsilon$ is the Euler vector field. If $X\in\cX^0(\Pi V)$,
$X\in\ker\pi$, then $X$ vanishes on $C^\infty(M)$. It is well known that there are vector fields $ {Z}_i\in\cX(M)$ and
smooth functions $f_i\in C^\infty(M)$ such that $\sum_i({Z}_i)(f_i)=1$ \cite{Grabowski}. As $X(f_i)=0$ by assumption, we can write
$$X=\left(\sum_iZ_i(f_i)\right)X=\sum_i\left([Z_i,f_iX]-[f_iZ_i,X]\right)$$
and it is clear that $\mathrm{supp}\,f_iX  \subset\mathrm{supp}\,X$. We can therefore assume the decomposition
(\ref{L}) with $\mathrm{supp}\,Y_i\subset\mathrm{supp}\,X$. Consider now a function $\psi\in C^\infty(M)$ with
the support in $U$ that takes value 1 on $\mathrm{supp}\,X$. We have
$$\sum_i[\psi X_i,Y_i]=\sum_i\left(\psi[X_i,Y_i]\pm(Y_i\psi)X_i\right)=\psi\sum_i[ X_i,Y_i]=\psi X=X\,,$$ since
$Y_i(\psi)=0$, as $\psi=1$ on the support of $Y_i$.
\end{proof}

Our next goal is to explain how the supermanifold
structure of $\cM$ is encoded in the Lie superalgebra structure of
$\cD^1(\cM)$ and, further, to describe all the automorphisms of
this superalgebra. More precisely, we will show that every Lie superalgebra isomorphism
$\cD^1 (\cM_1)\to\cD^1 (\cM_2)$, where
the dimension of manifolds is different from $0|1$, respects the filtration $\cD^0 (\cM_i)\subset
\cD^1 (\cM_i)$, $i=1,2$ and covers a supermanifold diffeomorphism $\cM_1\to\cM_2$.
We prove that
the group of automorphisms of the Lie superalgebra $\cD^1 (\cM)$ is a semi-direct product
of the subgroup induced by supermanifold diffeomorphisms and
another subgroup which consists of automorphisms determined by a 1-cocycle on the Lie superalgebra
of vector fields $\cX (\cM)$ with values in the algebra of functions on $\cM$.

\begin{rem} In contrast with even manifolds,
the Lie superalgebra structure on the Lie superalgebra of all
differential operators $\cD (\cM)$ for a general supermanifold $\cM$
does not recognize $\cD (\cM)$ as the Lie superalgebra of superdifferential operators on $\cM$: there are Lie superalgebra
isomorphisms $\cD (\cM_1)\to\cD (\cM_1)$ not respecting the canonical filtration of differential
operators $\cD^k (\cM)\subset \cD^{k+1}(\cM)$, $k\ge 0$. This is because, in the pure odd situation, any linear operator on superfunctions is a differential operator.
For instance, let $\cM$ be $\Pi V$ for
a finite-dimensional vector space $V$. Then $\cD (\cM)=\underline{\operatorname{End}}
(\Lambda^* (V))$ where $\Lambda^* V$ is viewed as a supervector space with the canonical
$\Z_2-$grading. Indeed, let $\xi^i$ be linear coordinates on $V$ counted as odd variables.
It is clear that $\cD (\cM)$ is generated as a associative superalgebra
by $\xi^i$ and $\pr_{\xi^i}$ subject to the Clifford relations
\beqn
\xi^i\pr_{\xi^j} +\pr_{\xi^j}\xi^i =\delta^i_j\,.
\eeq
Thus $\cD (\cM)$ is isomorphic as an associative superalgebra (correspondingly, as a Lie superalgebra)
to the Clifford algebra (correspondingly, the underlying Lie superalgebra) of $V\oplus V^*$
supplied with the canonical pairing. The latter coincides up to an isomorphism with the superalgebra of all linear endomorphisms of the spinor module $\Lambda^*V$. On the other hand, $V$ can be replaced with any maximal
isotropic subspace of $V\oplus V^*$. For instance, one can interchange $V$ and $V^*$; this will
produce a Lie superalgebra automorphism of $\cD (\cM)$ which apparently breaks the filtration.
Therefore the automorphism have not so nice geometrical description that is available for the Lie superalgebra of first-order differential operators.
\end{rem}

\section{Algebraic characterization of functions}

This section provides a Lie algebraic characterization of
superfunctions inside first-order superdifferential
operators. In what follows, $\cA$, $\cX$, $\cD^1$ denote
the algebras of superfunctions, supervector fields, first-order
superdifferential operators of a smooth supermanifold $\cM$ of
dimension $p|q$.\medskip
\begin{rem} In this paper we can assume that the odd
dimension $q$ of $\cM$ is at least $1$ -- otherwise we investigate
the already studied purely even situation \cite{Grabowski-Poncin}.
\end{rem}

\begin{theorem}\label{distinguishing_functions}
If $\operatorname{dim}\cM$ differs from $0|1$, the algebra
$\cA\subset\cD^1$ is the unique maximal super Lie ideal of $\cD^1$
consisting of $\operatorname{ad}$-nilpotent elements. In particular, any automorphism of the Lie superalgebra $\cD^1$ preserves $\cA$, thus induces an automorphism of the Lie superalgebra $\cX=\cD^1/\cA$. \end{theorem}

\begin{proof} Apparently, $\cA\subset\cD^1$ is an ideal made up by
$\operatorname{ad}$-nilpotent elements. It thus suffices to prove
that any another ideal with the same property is contained in
$\cA$.\medskip

First we identify $\cD^1/\cA$ as a Lie superalgebra with $\cX$ via
the canonical isomorphism $\cD^1/\cA\ni [D]\to D-D1\in\cX$.
Further, we choose any diffeomorphism of supermanifolds between
$\cM$ and $\Pi V$, where $V$ is a certain vector bundle $V\to M$
of rank $q$. This choice provides a non canonical isomorphism
between the superalgebras of functions, $\cA$ and $\zG(\wedge
V^*)$, which implements a $\Z_2$-compatible $\N$-grading on $\cA$.
The adjoint action of the corresponding derivation or Euler vector
field $\e\in\operatorname{Der}_0\cA=\cX_0$ supplies $\cX$ with a
$\Z_2$-compatible $\Z$-grading $\cX^k=\{X\in\cX:[\e,X]=kX\}$,
$k\in\Z$, $k\ge -1$ \cite{GKP1}.

Suppose now that $\cJ$ is a super Lie ideal of $\cD^1$ made up by
$\operatorname{ad}$-nilpotent elements and denote by $p:\cD^1\to
\cD^1/\cA\simeq\cX$ the canonical surjective Lie superalgebra
morphism. Then $\cI:=p(\cJ)$ is a super Lie ideal of $\cX$ whose
elements are $\operatorname{ad}$-nilpotent as well and that is
moreover $\Z$-graded with respect to $\e$. To explain the last
claim it suffices to prove that, if an element $X\in \cI\subset\cX$
decomposes as $X=\sum_{j=1}^s X_{k_j}$, $X_{k_j}\in\cX^{k_j}$,
$i>j\Rightarrow k_i>k_j$, then $X_{k_j}\in\cI$ for all $j$. As
$Y^m:=ad_\e^m X\in\cI$ for all $m$, we obtain
 \beqn
  \sum_{j=1}^s k_j^m X_{k_j}=Y^m\in\cI,\,\hspace{3mm} m=0,\ldots, s-1\,.
 \eeq
This linear system can be inverted, since the corresponding
Vandermonde matrix is nondegenerate, so that all $X_{kj}$ are
actually in $\cI$. Furthermore, it follows from the proof of
\cite[Proposition 2]{GKP1} that $\cI\cap \cX^0=\{0\}.$ Indeed, if
$X\in \cI\cap \cX^0$, then $X\in\cI_0$, which is a Lie ideal of
$\cX_0$ that is made up by ad-nilpotent elements. The mentioned
proof then implies that $X\in \cX^0\cap \oplus_{i>0}\cX^{2i},$ so
that $X=0.$

Let now $X$ be a nonvanishing element of $\cI$ and let
$X_k\in\cX^k\cap\cI$ be one of its nonvanishing $\Z$-homogeneous
terms. It is always possible to find $Y_1,\ldots, Y_r\in \cX$ such
that
 \beq\label{multibracket}
  \left[Y_1,\ldots,[Y_r, X_k]\ldots \right]\in \cX^0-\{0\}.\,
 \eeq
Indeed, it is easily seen that $X_k$ must be $0$, if the preceding
multibracket vanishes, for $X_k\in\cX^k,$ $k\ge 1$, (resp. for
$X_k\in\cX^{-1}$) and for all $Y_i\in\cX^{-1}$ (resp.
$Y_i\in\cX^1$) (a problem arises only if $\dim\cM=0|1$ and the
degree of $X$ is $-1$, as then $\cX^1=\{0\}$). Since the
multibracket (\ref{multibracket}) is again in $\cI$, we get a
contradiction with the fact that $\cI\cap \cX^0=\{0\}$. Therefore
$p(\cJ)=\cI=\{0\}$, so that eventually $\cJ\subset\cA$.
\end{proof}
\begin{rem}
The assumption that the dimension of $\cM$ differs from $0|1$ is crucial. Indeed,
for a manifold of dimension $0|1$ there are automorphisms of the Lie algebra $\cD^1$ not preserving $\cA$, For, let $\xi$ be an odd coordinate on $\cM$. Then, $\cD^1$ is spanned by superfunctions $1,\xi$ and supervector fields $\p_\xi,\xi\p_{\xi}$. It is easy to see that the liner map on $\cD^1$ for which
$$\xi \p_\xi\mapsto -\xi \p_\xi\,,\quad \p_\xi\leftrightarrow\xi\,,\quad 1\mapsto 1$$
is an automorphism.
\end{rem}

\section{Reduction of the automorphism-problem}

In this section we fix a supermanifold $\cM$ of dimension different from $0|1$ and we reduce the quest for the automorphisms of the Lie algebra of first-order superdifferential operators to the
computation of the even 1-cocycles of the Lie algebra of
supervector fields valued in the associative algebra of superfunctions.\medskip

Clearly, any diffeomorphism $\zvf$ of the supermanifold $\cM$
induces an automorphism of the Lie superalgebra of first-order
differential operators of $\cM.$ Indeed, any
$\zvf\in\mathrm{Diff}(\cM)$ defines in particular an associative
superalgebra automorphism $\zvf^*:{\cal A}\to{\cal A}$ and the
induced Lie superalgebra automorphism is given by
$$\zvf_*:{\cal
D}^1\ni D\mapsto \zvf^{*{-1}}\circ D\circ \zvf^*\in{\cal
D}^1\,.$$

In the following we use the canonical splitting $\cD^1=\cA\oplus\cX$ and we denote the projection $\cD^1\ni D\mapsto
D-D1\in \cX$ by $\rho$. Using the fact that, for $D\in{\cal D}^1$
and $f\in{\cal A}$, we have $[D,f]=Df-(D1)f$, we immediately check
that $\zr$ is a representation by derivations of the Lie
superalgebra ${\cal D}^1$ on the supervector space $\cal A$.

\begin{proposition}\label{description} Any Lie superalgebra automorphism $\phi\colon \cD^1\to
\cD^1$ splits into a product $\varphi_*\circ \phi_c$, where
$\varphi_*$ is the automorphism induced by a diffeomorphism
$\varphi\colon\cM\to\cM$, and where the automorphism $\zf_c$ has
the form $\phi_c=\mathrm{id}+c,$ with $c$ being an even 1-cocycle of
$\cD^1$ represented upon $\cA$ by $\zr$.
\end{proposition}
\begin{proof} It follows from Theorem
\ref{distinguishing_functions} that any automorphism $\phi\colon
\cD^1\to \cD^1$ preserves $\cA$. Thus $\phi$ induces a Lie
superalgebra automorphism \bee\tilde\phi:{\cal X}\ni X\mapsto
\phi(X)-\phi(X)1\in {\cal X}.\label{AutVfInd0}\eee 

All automorphisms of $\cal X$ are classified in
\cite{GKP1}, see Theorem 2 and Proposition 5. It follows from the
proof of Proposition 5 that they are implemented by a
diffeomorphism of $\cM$, unless $\dim\cM=1|1$ or
$0|2$. In the latter case there exist additional
automorphisms. Let us assume that $\dim\cM=1|1$ or
$0|2$. Then the Euler vector field $\e\in{\cal X}_0$ is
well-defined up to a sign, such that any automorphism of $\cal X$,
which is not coming from a diffeomorphism of $\cM$, interchanges
$\e$ with $-\e$.
Now we can make a canonical choice of the Euler vector field $\e$ such
that $\e$ acts on $\cA_1$ as the identity. Apparently, any automorphism
of $\cal X$, which is induced by an automorphism of $\cD^1$, must preserve the canonically
chosen $\e$. Thus, in both exceptional (low-dimensional) cases an automorphism
of $\cal X$ which
exchanges $\e$ and $-\e$ does not admit an extension to first-order
differential operators. Therefore, each $\tilde\phi$ is always induced by
a diffeomorphism $\varphi\colon\cM\to\cM$.

Let us take a diffeomorphism $\varphi$ of $\cM$, such that
$$\tilde\phi:\cX\ni X\mapsto \zvf^{*{-1}}\circ X\circ \zvf^*\in\cX.$$ The automorphism $\tilde
\psi=\tilde\zvf_*^{-1}\circ\tilde\zf$ of $\cal X$, induced by the
automorphism $\psi:=\zvf_*^{-1}\circ\zf$ of ${\cal D}^1$, is the
identity map. Hence, for $X\in\cX$, we have $\psi(X)=X+\psi(X)1$
and, for $D=f+X\in{\cal D}^1$, we thus get $\psi
(D)=f+X+\psi(X)1+\psi(f)-f=D+c(D)$, where $c:{\cal D}^1\to \cA$ is
an even linear map. The injectivity and surjectivity of $\psi$ are
equivalent to the corresponding property of
$(\mathrm{id}+c)|_{\cA}$. To prove that the surjectivity of this
restriction implies that of $\psi$, it suffices, if $\zD=g+Y\in
\cD^1$, to set $g-c(Y)=f+c(f)$, for some $f\in\cA$, and to observe
that $\psi(f+Y)=\zD$. Eventually, the Lie algebra homomorphism
property of $\psi=\mathrm{id}+c$ is clearly equivalent with the
1-cocycle condition of $(\cD^1,\zr)$ for $c$. \end{proof}

\begin{theorem}\label{Automorphism2} The automorphism $\phi_c$ can be uniquely written in the form
$$\phi_c(f+X)=(\zk f+\gamma(X))+X\,,$$ where $\zk$ is a non-zero constant and
$\zg:\cX\rightarrow\cA$ is an even 1-cocycle of the Lie superalgebra $\cX$ canonically represented on
$\cA$.\end{theorem}

\begin{proof} Note first that, in view of the preceding proof, any even $1$-cocycle
$c$ of the Lie superalgebra $\cD^1$ represented by $\zr$ on $\cA$,
such that $(\mathrm{id}+c)|_{\cA}$ is bijective, defines an
automorphism $\zf_c=\mathrm{id}+c$. The 1-cocycle condition
obviously splits into the intertwining condition \bee
c(Xf)=X(c(f)),\;\forall X\in\cX, f\in\cA,\eee and the 1-cocycle
condition of the Lie superalgebra $\cX$ canonically represented
upon $\cA$. We prove below in Lemma \ref{Intertwin1} that the intertwining condition means that
$c\,|_{\cA}=\zl\cdot\mathrm{id}$, $\zl\neq -1$, where the exclusion of the value $-1$ is due to the condition
that $(\mathrm{id}+c)|_{\cA}$ be a bijection. Thus, $$\phi_c(f+X)=f+X+c(f)+c(X)=((1+\lambda)f+c(X))+X$$ and it is easy to see that $\zg:\cX\rightarrow\cA$, $\zg=c_{|\cX}$, is a 1-cocycle with coefficients in the canonical representation.
  \end{proof}

\begin{lemma}\label{Intertwin1} Any even linear map $c:\cA\to \cA$ that verifies the intertwining condition $$c(Xf)=X(c(f))$$ for all $X\in\cX$ and $f\in\cA$ is of the form $c(f)=\lambda f$, $\lambda\in\R$.\end{lemma}

\begin{proof} Like in the proof of Theorem \ref{local}, we can assume that $\cM=\zP V$, where $V$ is a vector bundle over $M$, and consider vector fields $Z_i\in\cX(\zP V)$ and functions $f_i\in C^{\infty}(M)$, such that $\sum_i Z_i(f_i)=1$.
Since $$f=f\cdot 1=\sum_i(fZ_i)(f_i)\,,$$ we see that $c(f)=f\sum_i Z_i(c(f_i))$. Let us put $\zl=\sum_i Z_i(c(f_i))\in\cA$. As $\zl=c(1)$ and for each $X\in\cX$ we have $X(\zl)=X(c(1))=c(X(1))=0$, the superfunction $\zl$ is actually a constant.
\end{proof}
\begin{rem} The additional automorphisms of the Lie superalgebra $\cX(\cM)$ in the exceptional cases can be described as follows.
In the first case,  when $\dim\cM=1|1$, the supermanifold
$\cM$ is isomorphic to $\Pi L$ for a real line bundle $L\to M$
and $\cA_1\simeq\Gamma (L^*)$. Taking into account that the structure group
$L$ can be reduced to $O(1)=\Z_2$, one can always choose a trivialization
$\sigma_0$ of $L^{\otimes 2}$ and a flat connection $\na$ on $L$
(and thus a flat connection on $L^{\otimes 2}$, denoted by the same letter) such that
$\na (\sigma_0)=0$. On the other hand, $M$ is either $\R$
or $S^1$, so $M$ is always orientable. Let $\mu$ be a volume form on
$M$, then $\sigma_0\otimes\mu$ determines a bundle isomorphism $L^* \otimes TM\to L$,
which gives rise to a bundle automorphism $\chi$ of $L\oplus L^*\otimes TM$,
interchanging $L$ and $L^*\otimes TM$ such that $\chi^2 =\Id$. Taking into
account that $\cX_1$ is isomorphic to $\Gamma (L\oplus L^*\otimes TM)$ as a vector space,
we obtain an invertible linear map $\cX_1\to\cX_1$.
The choice of a flat connection $\na$ determines
a Lie algebra isomorphism $\cX_0\simeq D^1 (M)$, where $D^1 (M)$ is the Lie algebra of first-order differential operators on $M$. The commutator relations in $\cX$ are given by the following
formulas:
\beqn
[v+f, s]&=&\na_v (s)-fs\,\\\nonumber
[v+f, \theta\otimes v'] &=& \na_v (\theta)\otimes v' +\theta\otimes
[v,v']+f \theta\otimes v'\,,\\\nonumber [s, \theta\otimes v']&=&\langle\theta , s\rangle v' +\langle\theta,\na_{v'}s\rangle\,,
\eeq
where $v, v'\in\Gamma (TM)$, $s\in\Gamma (L)$, $\theta \in \Gamma (L^*)$, and $f\in C^\infty (M)$.
Apparently, $\cX_1$ is a faithful ${\cal X}_0-$module. The linear
invertible map $\chi$ determines a Lie algebra automorphism of  ${\cal X}_0$
of the form $v+f\mto v + div_\mu (v)-f$, where $div_\mu (v)=L_v (\mu)\mu^{-1}$.
Using appropriate local
supercoordinates $(t,\xi)$,
such that $\mu =\md t$, $\na (\xi)=0$, and $\sigma_0 =\xi^{\otimes (- 2)}$,
 we can write the whole transformation in the following form (here
$h$ and $f$ are arbitrary local smooth functions):
\beqn
h(t)\pr_\xi &\mto& h(t)\xi\pr_t\,,\\ \nonumber
h(t)\xi\pr_t &\mto& h(t)\pr_\xi\,,\\ \nonumber
h(t)\pr_t + f(t)\xi\pr_\xi &\mto& h(t)\pr_t + (\pr_t h(t) -f(t))\xi\pr_\xi\,.
\eeq
It is easy to verify that such a transformation is a Lie superalgebra automorphism.
In fact, the group of supermanifold diffeomorphisms of $\cM$ acts
freely and transitively on the set of automorphisms which
interchanges $\e$ and $-\e$, so a particular choice of $\chi$
is nothing but the choice of "an origin".

In the second case, when $\dim\cM=0|2$, the supermanifold $\cM$ is isomorphic to $\Pi V$ for a 2-dimensional
vector space $V$. The even part of $\cal X$ is naturally isomorphic to $\gl (V)$,
while the odd part is isomorphic to $V\oplus \Lambda^2 V^*\otimes V$. Let us
fix a constant volume form $c\in \Lambda^2 V^*$. The Lie algebra $\gl (V)$ is a
direct sum of $\sll (V)$ and the center of $\gl (V)$ spanned by $\Id$. The subalgebra $\sll (V)$
preserves $c$, which makes $V$ and $\Lambda^2 V^*\otimes V$ into isomorphic $\sll (V)$-modules. Let us combine the isomorphism $V\stackrel{c\otimes }{\longrightarrow}
\Lambda^2 V^*\otimes V$ with the Lie algebra isomorphism of $\gl (V)$ which preserves
all elements of
$\sll (V)$ and interchanges $\Id$ and $-\Id$ (such an automorphism is a composition of a conjugation
and the opposite to a transposition). One can easily check that the obtained
linear map $\cX\to\cX$ is a Lie superalgebra automorphism.
\end{rem}

\section{Cohomology of supervector fields represented on functions}

Let $\zg:\X(\cM)\rightarrow\cA(\cM)$ be an even 1-cocycle, so
that, for any $X,Y\in\cX(\cM)$,
\beq\label{cocycle}\zg([X,Y])=X(\zg(Y))-(-1)^{\vert
X\vert\vert Y\vert}Y(\zg(X))\,.\eeq

\begin{prop} Every even 1-cocycle of the Lie superalgebra $\cX(\cM)$ represented upon $\cA(\cM)$ is a local operator. \end{prop}

\begin{proof} Let $X\in\X(\cM)$ such that $X|_U=0$,
$U\subset M$ open, and let $x_0\in U$. According to Theorem \ref{local}, the vector field $X$ reads
$X=\sum_i[X_i,Y_i],$ for some $X_i\in\X(\cM)$ and some $Y_i\in\X_0(\cM)$ which are 0 in a neighbourhood of
$x_0$ as well. Hence, in view of the cocycle condition,
$$\zg(X)=\sum_i\left( X_i(\zg( Y_i))-(-1)^{\vert
X_i\vert\vert Y_i\vert} Y_i(\zg( X_i))\right)\,,$$
$\zg(X)=0$ in a neighbourhood of $x_0$.
\end{proof}

The next result gives the local form of the even 1-cocycles of
$\cX(\cM)$ with coefficients in $\cA(\cM)$. We will denote by $\zW^1(U)$, $U\subset M$ open, the super $\cA(U)$-module of superdifferential 1-forms over $U$
and $\zW^1_0(U)$ will refer to its even part.

\begin{theorem}\label{Cohomology0Loc} Let $(U,u=(u^1,\ldots,u^{p+q}))$ be any coordinate chart of $\cM$. The restriction to $U$ of any even 1-cocycle
$\zg:\X(\cM)\rightarrow\cA(\cM)$ is of the form \beq\label{D}
\zg|_U\left(\sum_k\partial_k\cdot
g^k\right)=\sum_k\left(a\,\partial_k g^k+\omega_k\,g^k\right)\,,
\eeq where $\p_k=\p_{u^k}$, $g^k\in\cA(U)$, $a\in\R$, and
$\omega=\sum_k du^k\,\omega_k$ is a closed even 1-form.
\end{theorem}

\begin{rem} As in this text we use the Deligne sign convention for the wedge product, the super de Rham operator $d$ has parity 0.\end{rem}

\begin{proof} Due to locality, the restriction
$\zg|_U$ is an even 1-cocycle of $\cX(U)$ valued in $\cA(U)$. In
the following we often omit the restriction to $U$ and write
simply $\zg$, $\cX$, $\cA$, etc. We will of course show that $\zg$ is
a differential operator.\medskip

When looking at the cocycle condition for $X=u^i\p_j$ and
$Y=g\p_k$, where $g\in \cA$, we are naturally led to introduce the
map
$$\zg_k:\cA\ni g\mapsto
(-1)^{|u^k||g|}\zg(g\partial_k)=\zg(\partial_k\cdot g)\in\cA\,.$$
The cocycle equation then means that the map
\beq\label{e1}g\mapsto
\zg_k(u^i\pr_jg)-\delta^i_k\zg_jg-(-1)^{|u^k|(|u^i|+|u^j|)}u^i\pr_j(\zg_kg)
\eeq is a differential operator of order 0 and parity
$|u^i|+|u^j|+|u^k|$. Similarly, taking $X=\pr_j$ and $Y=g\pr_k$,
we get that the map \beq\label{e2} g\mapsto
\zg_k(\pr_jg)-(-1)^{|u^k||u^j|}\pr_j(\zg_kg) \eeq is a
differential operator of order 0 and parity $|u^j|+|u^k|$. Thus,
subtracting from operator (\ref{e1}) the operator (\ref{e2})
multiplied from the left by $(-1)^{|u^k||u^i|}u^i$, we obtain that
\beq\label{e3}T_{k,i,j}=[\zg_k,u^i]\circ\pr_j-\delta^i_k\zg_j\,
\eeq is a differential operator of order 0 and parity
$|u^i|+|u^j|+|u^k|$.\medskip

If $i\ne k$, then $T_{k,i,j}$ reduces to
$[\zg_k,u^i]\circ\partial_j$. The latter is of order 0 and
vanishes on constants, so $[\zg_k,u^i]\circ\partial_j=0$, for all
$i\ne k$ and all $j$. This in turn implies that
\beq\label{o1}[\zg_k,u^i]=0\quad\text{for}\quad i\ne k\,. \eeq
Indeed, if there exists an even coordinate $u^j$, we can integrate
classical functions with respect to $u^j$ and the result is
obvious. In the pure odd case, we have in particular
$[\zg_k,u^i]\circ\partial_i=0\,,$ $i\neq k$, so that
\beq\label{op3}0=[\zg_k,u^i](\partial_i(u^ig))=[\zg_k,u^i](g)+(-1)^{\vert
u^i\vert}[\zg_k,u^i](u^i\partial_ig),\eeq for all $g\in\cA$. But,
as easily checked, for odd $u^i$, the supercommutator
$[\zg_k,u^i]$ supercommutes with the multiplication by $u^i$, so
that
$$[\zg_k,u^i]\circ(u^i\partial_i)=(-1)^{|u^i|(|u^k|+|u^i|)}u^i[\zg_k,u^i]\circ\partial_i=0,$$
and, according to Equation (\ref{op3}), $[\zg_k,u^i]=0$, which
completes the proof of Equation (\ref{o1}).\medskip

For $i=k$, the differential operator $T_{i,i,j}$ reads
$$T_{i,i,j}=[\zg_i,u^i]\circ\partial_j-\zg_j$$ and, since it
is of order 0, we have $[T_{i,i,j},u^j]=0$, i.e.,
\beq\label{order1}[\zg_i,u^i]+(-1)^{|u^j|}[[\zg_i,u^i],u^j]\circ
\pr_j-[\zg_j,u^j]=0\,. \eeq As for $i\ne j$, the Jacobi identity
and Equation (\ref{o1}) entail
$$[[\zg_i,u^i],u^j]=(-1)^{|u^i||u^j|}[[\zg_i,u^j],u^i]=0\,,
$$
we get \beq\label{1order} [\zg_i,u^i]=[\zg_j,u^j]\,. \eeq Choosing $i=j$ in Equation (\ref{order1}), we obtain
$[[\zg_i,u^i],u^i]\circ \pr_i=0,$ which implies \beq\label{1order1}[[\zg_i,u^i],u^i]=0\,,
\eeq if $u^i$ is even. However, it suffices to develop the {\small
LHS} of Equation (\ref{1order1}) to conclude that the claim holds
true for odd $u^i$ as well. Therefore, when taking into account
Equation (\ref{o1}), we finally get
\beq\label{1order2}[[\zg_i,u^j],u^k]=0\,, \eeq for all
$i,j,k$. The latter equation means that $\zg_k$ are first-order
differential operators (see Remark \ref{order}).

According to Equation (\ref{o1}), $\zg_k$ commute with
multiplication by $u^i$, $i\ne k$, so they are of the form
$\zg_k=a_k\pr_k+\omega_k$, $a_k,\zw_k\in\cA$. Since, for any first-order operator and any function, we have $[D,f]=Df-D1\cdot f$,
Equation (\ref{1order}) shows that $a_i=a_j=:a$. Finally,
\beq\label{form} \zg\left(\sum_k\pr_k\cdot
g^k\right)=\zg\left(\sum_k (-1)^{|u^k||g^k|}
g^k\pr_k\right)=\sum_k\left(a\,\pr_kg^k+\omega_k\,g^k\right)\,.
\eeq Note now that $\vert a\vert=0$ and $\vert \omega_k\vert=\vert
u^k\vert$, since $\zg_k$ has parity $|u^k|$. Starting from
$$[f\partial_i,g\partial_j]=f\p_ig\,\partial_j
-(-1)^{(\vert u^i\vert+\vert f\vert)(\vert u^j\vert+\vert g\vert)}
g\p_jf\,\partial_i\,,$$ we straightforwardly see that the
corresponding cocycle condition provides, after simplification, an
identically vanishing bidifferential operator in $f$ and $g$. When
writing that its coefficients vanish, we get \beq
\p_ia&=&0\,,\\\p_i \omega_j-(-1)^{\vert u^i\vert\vert
u^j\vert}\p_j \omega_i&=&0\,, \eeq for all $i,j$. We thus conclude
that the differential operator $\zg$ defined by (\ref{form}) is a
1-cocycle if and only if $a\in\R$ and the even superdifferential
1-form $\omega=\sum_k du^k\omega_k$ is closed.
\end{proof}

It is easily checked that $\zg=\zg|_U$ defined by
Equation (\ref{form}) has the following property with respect to
the right module structure of $\cX=\cX(U)$ over $\cA=\cA(U)$:
\beq\label{gdiv} \zg(X\cdot f)=\zg(X)\cdot f+a\,Xf\;.\eeq

\begin{prop} Any even 1-cocycle $\zg:\cX(\cM)\to \cA(\cM)$ verifies Equation
(\ref{gdiv}), for all $X\in\cX(\cM)$, all $f\in\cA(\cM)$, and some $a\in\R$.\end{prop}

\begin{proof} Indeed, the restriction $\zg|_U$ to any chart domain $U$ is of the
form (\ref{form}) and thus satisfies the local equation
(\ref{gdiv}) for some $a_U\in\R$. It follows that, for any global
$X$ and $f$, we have
$$a_U(Xf)|_{U\cap V}=(\zg(X\cdot f)-\zg(X)\cdot f)|_{U\cap
V}=a_V(Xf)|_{U\cap V},$$ where $V$ denotes a chart domain that
intersects $U$. Since the base manifold $M$ of $\cM$ is connected,
all $a_U$ coincide, which proves the claim.  \end{proof}

Some authors, see e.g. \cite{KSM}, define a divergence operator in $\cM$ as an operator $\zg:\X(\cM)\rt\cA(\cM)$
that satisfies only (\ref{gdiv}) with $a=1$. In this paper, we assume
also the cocycle condition, which can be understood
as a vanishing curvature condition for $\zg$.

\begin{deff} A {\em divergence} in a supermanifold $\cM$ (or a {\em superdivergence}) is a 1-cocycle of $\cX(\cM)$ valued in $\cA(\cM)$ that satisfies Equation (\ref{gdiv}) with $a=1$.
Any even 1-cocycle will be called a {\em generalized divergence}.\end{deff}

\begin{rem} We prove in the last section that in any smooth supermanifold there exists a divergence
operator $\zg_0$.\end{rem}

We are now able to describe the even part of the first cohomology group
of the Lie superalgebra $\cX(\cM)$ valued in $\cA(\cM)$.

\begin{theorem}\label{Automorphism3} Let $\cM$ be a supermanifold with a fixed divergence $\zg_0$. Any even 1-cocycle $\zg:\X(\cM)\rt\cA(\cM)$, i.e.
any generalized divergence, can be
uniquely written as
$$\zg=a\,\zg_0+i_{\omega}\,,$$
where $a\in\R$ and $\omega$ is a closed even 1-form on $\cM$. The cocycle $\zg$ is a coboundary if and only if $a=0$ and
$\omega=\text{d}f$, $f\in\cA_0(\cM)$. In other words, $$H^1_0(\cX(\cM),\cA(\cM))=\R\,\zg_0\oplus H^1_{\mathrm{DR},0}(\cM),$$ where $H^1_{\mathrm{DR},0}(\cM)$ is
the even part of the first super de Rham cohomology group of $\cM$.
\end{theorem}

\begin{proof} If $\zg$ is a generalized divergence, it verifies Equation (\ref{gdiv}) for some $a\in\R$. Therefore, the difference $\zg-a\zg_0$ is a right
$\cA(\cM)$-module morphism from $\cX(\cM)$ to $\cA(\cM)$, so an even super 1-form $\zw\in\zW^1_0(\cal M)$. The cocycle equation for $\zw=\zg-a\zg_0$ and
Cartan's formula for the super de Rham operator $d$ show that $d\zw=0.$ Furthermore, $a\zg_0+\zw=\zg=df$, $f\in\cA_0(\cM)$, if and only if $a=0$ and $\zw=df.$  \end{proof}

It follows from Proposition \ref{description} and Theorems
\ref{Automorphism2} and \ref{Automorphism3} that for
supermanifolds we get an analog of Theorem 8 from
\cite{Grabowski-Poncin}.

\begin{theorem} Let $\cal M$ be a smooth supermanifold of dimension different from $0|1$ and let $\zg_0$ be a fixed divergence in $\cM$. Then, an even linear map  $\phi:\cD^1(\cM)\to\cD^1(\cM)$ is an automorphism
of the Lie superalgebra $\cD^1(\cM)=\cA(\cal M)\oplus\cX(\cM)$ of first-order superdifferential operators of $\cM$ if and only if it can be written in the form
\begin{equation}
\phi(f+X)=\zvf_*(X)+(\zvf^{-1})^*(\kappa f+a\,\zg_0(X) +i_{\omega}(X)),
\end{equation}
where $\zvf$ is a diffeomorphism of $\cM$, $a,\kappa\in\R$, $\kappa\ne 0$, and $\omega$ is a closed even 1-form. All the objects $\zvf,a,\kappa,$ and $\omega$ are uniquely determined by $\phi$.
\label{d1}
\end{theorem}

\section{Existence of superdivergences}
Our aim in this section is to prove the existence of a divergence on each supermani\-fold. This will be done similarly to the proof of existence of a divergence on any standard (even) manifold with the use of a nowhere-vanishing 1-density understood as a class of volume forms up to a sign. The sheaf of top forms $\zW^{\mathrm{top}}$ of a classical
differential manifold has to be replaced, in the case of a supermanifold ${\cM}=(M,\cA)$, by the {\it Berezinian sheaf}
$\mathrm{Ber}=\mathrm{Ber}(\cM)$ whose nowhere-vanishing sections are {\it Berezinian volumes}. A homogeneous Berezinian volume $s$ defines a divergence $\gamma_s$ of a homogeneous vector field $X$ by the formula (see \cite{KSM})
$$ \cL_Xs=(-1)^{|X||s|}s\cdot\gamma_s(X)\,.$$
Here $\cL_Xs$ is the Lie derivative of the Berezinian volume understood as a differential operator on $\cA$ defined by $\cL_Xs=-(-1)^{|X||s|}s\circ X$.\medskip

Let us fix an atlas of supercharts of $\cal M$.
Over a domain $U$ with supercoordinates $u=(x^1,\ldots,x^p,\xi^1,\ldots,\xi^q)$,
the (right) $\cA(U)$-module $\mathrm{Ber}(U)$ is given by
$$\mathrm{Ber}(U)=\zG(U,\mathrm{Ber}):=d^{\,p|q}u\,\cA(U)\simeq
\left(dx^1\wedge\ldots\wedge dx^p\otimes
\p_{\xi^q}\circ\ldots\circ\p_{\xi^1}\right)\,\cA(U)\,.$$
Since every supermanifold is diffeomorphic to $\Pi V$ for a vector bundle $V$ of rank $q$ over $M$, we can assume that $\cM=\Pi V$, so that $\cA$ is a superalgebra isomorphic with the Grassmann algebra of sections of $\Lambda^\bullet V^*=\bigoplus_{k=0}^\infty\Lambda^kV^*$ and local odd coordinates $(\xi^1,\ldots,\xi^q)$ are represented by sections of $V^*$. One can easily see that with any section $v$ of the line bundle $\Vol(V)=\Lambda^pT^*M\otimes_M\Lambda^q V$ one can associate a Berezinian volume $s_v$. In other words, we have an embedding $S$ of the space of sections of $\Vol(V)$ into the space of Berezinian densities. In local coordinates this embedding reads
$$dx^1\wedge\ldots\wedge dx^p\otimes e_q\wedge\ldots\wedge e_1\ \mapsto \ dx^1\wedge\ldots\wedge dx^p\otimes
\p_{\xi^q}\circ\ldots\circ\p_{\xi^1}\,,$$
where $e_1,\ldots,e_q$ is a basis of local sections of $V$ and $\p_{\xi_i}$ is the derivation of $\cA$ being the contraction of a section of $\Lambda^\bullet V^*$ with $e_i$. A linear change of coordinates in the vector bundle $V$, say
$(x,\xi)\mapsto(y(x),\eta)$ with $\eta_j=\sum_i a^i_j(x)\xi_i$, results in multiplying both sides by
$$\det\left(\frac{\p y^a}{\p x^b}(x)\right)\det\left(a^i_j(x)\right)^{-1}\,,$$
so the embedding is well defined globally. Since the line bundles are classified by $H^1(M,\Z_2)$, we can find an open covering $(U_\alpha)$ on $M$ and local nowhere-vanishing sections $v_\alpha$ of $\Vol(V)_{|U_\alpha}$ such that $v_\alpha=\pm v_{\alpha'}$ on $U_\alpha\bigcap U_{\alpha'}$. Hence, $S(v_\alpha)=\pm S(v_{\alpha'})$ over $U_\alpha\bigcap U_{\alpha'}$. But the divergence $\gamma_s$ does not depend on the sign of $s$, so the collection $\sigma$ of local nowhere-vanishing Berezinian volumes $S(v_\alpha)$ gives rise to a super(divergence) $\gamma_\sigma$ on $\cM$. Note that the collection of local nowhere-vanishing Berezinian volumes $S(v_\alpha)$ can be viewed as a nowhere-vanishing section of the {\it $1$-density sheaf} ${\frak D}_1=\Ber\otimes\mathrm{or}(M)$ of ${\cal M}$, defined as
the Berezinian sheaf twisted by the orientation sheaf of the body.
In the literature, the sections of ${\frak D}_1$ are sometimes
referred to as {\it nonoriented Berezinian sections}. If $\sigma'=f\sigma$,  with an even nowhere-vanishing function $f$, is another such section, then $\zg_{\sigma'}(X)=\zg_\sigma(X)+X(\ln
 f)$. Thus we get the following.

 \begin{theorem} Let $\cM=(M,\cA)$ be a smooth supermanifold. Then, there is a nondegenerate global section $\sigma$ of ${\frak D}_1$ with even coefficients. Moreover, for any such $\sigma$,
 $$\zg_\zs:\cX\ni X\mto (\cL_X \sigma)\sigma^{-1}\in\cA$$ is a divergence of $\cM$. Any other such section of ${\frak D}_1$
 implements a divergence in the same cohomology class.
 \end{theorem}

\vskip1cm
\noindent Janusz GRABOWSKI\\ Polish Academy of Sciences\\ Institute of
Mathematics\\ \'Sniadeckich 8, P.O. Box 21\\ 00-956 Warsaw,
Poland\\Email: jagrab@impan.pl \\

\noindent Alexei KOTOV\\ University of Luxembourg\\
Facult\'e des Sciences, de la Technologie et de la Communication\\
6, rue Richard Coudenhove-Kalergi\\
L-1359 Luxembourg City, Grand-Duchy of Luxembourg\\
Email: alexei.kotov@uni.lu\\

\noindent Norbert PONCIN\\University of Luxembourg\\
Facult\'e des Sciences, de la Technologie et de la Communication\\
6, rue Richard Coudenhove-Kalergi\\
L-1359 Luxembourg, Grand-Duchy of Luxembourg\\Email: norbert.poncin@uni.lu

\end{document}